\documentclass[a4paper, 10pt]{article}
\usepackage{amscd, amssymb}
\usepackage[mathscr]{eucal}
\usepackage[all]{xy}
\usepackage{mathrsfs}
\usepackage{amsfonts}
\usepackage{amsmath}
\usepackage{amsthm}
\usepackage{latexsym}
\usepackage{graphicx}
\newcommand{\Nv}{\vec{N}}
\newcommand{\N}{\scriptscriptstyle{N}}

\newcommand{\V}{\mathbf{v}}
\newcommand{\rk}{\mathrm{Rank}}
\newcommand{\R}{\mathbf{r}}
\newcommand{\T}{\mathbf{t}}
\newcommand{\Gg}{\mathbf{\Gamma}_N}
\newcommand{\h}{H_{1,k,c}(\Gg)}
\newcommand{\di}{\mathrm{dim}}
\newcommand{\Hom}{\mathrm{Hom}}
\newcommand{\C}{\mathbf{C}}
\newcommand{\real}{\mathbf{R}}
\newcommand{\s}{\sigma}

\newcommand{\g}{\gamma}
\newcommand{\G}{\Gamma}
\newcommand{\pf}{\textsl{Proof.}}
\newcommand{\epf}{\begin{flushright} $\Box$ \end{flushright}}
\newcommand{\rem}{\textbf{Remark.}}

\newtheorem{defi}{Definition}[section]
\newtheorem{thm}[defi]{Theorem}
\newtheorem{lem}[defi]{Lemma}
\newtheorem{prop}[defi]{Proposition}

\newcommand{\M}{\scriptscriptstyle {M}}
\newcommand{\e}{\epsilon}

\usepackage{latexsym}

\begin{document}

\title{On some finite dimensional representations of symplectic
reflection algebras associated to wreath products}

\author{Silvia Montarani}

\maketitle

\section{Introduction}\label{intro}

 Let $\G$ be a finite subgroup of $SL(2,\C)$, $k\in\C$ a complex number and  $c=\{c_{\g}\}$ a  complex valued class  function on $\G-\{e\}$, where $e$ is the identity element of $\G$. In this note we study some finite dimensional representations of the wreath product symplectic reflection algebra $\h$  attached to the group $\Gg:=S_N\ltimes \G^N$  and to the parameters $(k,c)$ (see \cite{EG},\cite{GG}).

We recall that in the rank $1$ case there is no parameter $k$ and finite
dimensional representations of the wreath product algebra $H_{1,c}(\G)$ have
been classified in \cite{CBH}. 
For higher rank,  when  $k=0$, we have that 
$H_{1,0,c}(\Gg)=S_N\sharp H_{1,c}(\Gamma)^{\otimes N}$.
 So, in this particular case, the finite dimensional representations of $H_{1,k,c}(\Gg)$ can also be recovered from \cite{CBH}.

Thus we restrict ourselves to values of $(k,c)$ with  $k\neq 0$ and obtain a classification of the irreducible finite dimensional representations of $\h$ in which the generators $x_i$ and $y_i$  (see Definition \ref{sympl}) act by $0$ for any $i=1,\dots,N$. We show that this is equivalent to classifying all the irreducible representations of $\Gg$ that can be  extended to representations of $\h$ for such values of the parameter $(k,c)$ .

\section{Preliminaries}\label{preli}

\subsection{Symplectic reflection algebras of wreath product type}\label{symal}

In this first section we will give the definition of the wreath product symplectic reflection algebra $\h$ using an explicit presentation by generators and relations due to  \cite{GG}. We start by recalling the following well known definition. 
Let $\G$ be a finite group and  $\G^N$ be the direct product of $N$ copies of $\G$. Denote by $S_N$ the symmetric group of degree $N$.

\begin{defi}\label{wreath}
The wreath product  $\Gg:=S_N\ltimes \G^N$ is the semi-direct product of $S_N$ and $\G^N$, where $\G^N$ is normal, and the action of $S_N$ on $\G^N$ by conjugation is the natural one in which $S_N$ permutes the direct factors of $\G^N$. 
\end{defi}

Let now $L$ be a $2$-dimensional complex vector space with a symplectic form
$\omega_L$, and consider the space $V=L^{\oplus N}$, endowed with the induced symplectic form $\omega_V={\omega_L}^{\oplus N}$.  Chosen a symplectic basis $x$, $y$ for $L$ we will denote by $x_i$, $y_i$ the corresponding vectors in the $i$th $L$-factor of $V$.
 
Let $\Gamma$ be a finite subgroup of $Sp(L)$ ( note that $Sp(L)=SL(2,\C)$ thanks to the choice of a symplectic basis), then $\G^N$ acts symplectically on $V$. Consider the action of  $S_N$   on $V$ that permutes the factors in the direct sum. The group $\Gg:=S_N\ltimes\Gamma^N \subset Sp(V)$ acts naturally on $V$.  We will denote by $e$  the identity element of $\G$  and we will write  $\gamma_i$ for an element $\gamma\in \Gamma$, seen as the element $(e,\dots, \stackrel{i}{\g},\dots, e)$ in the $i$th factor of $\G^N\subset \Gg$.

Let $(k,c)$ be a pair made up by a complex number $k$ and a complex valued class function  $c$ on $\G-\{e\}$. We will write  $c_{\g}$ for $c(\g)$. 

We denote by $TV$ the tensor algebra of $V$ and by $\Gg\sharp TV$ the smash product of $TV$ with $\Gg$.

\begin{defi}\label{sympl}The algebra $\h$ is the quotient
of $\Gg\sharp TV$ by the following relations:
\begin{itemize}
\item[{\emph{(R'1)}}] 
For any $i\in [1,N]$: $$[x_{i}, y_{i}]
=  1+ \frac{k}{2} \sum_{j\neq i}\sum_{\g\in\Gamma}
s_{ij}\g_{i}\g_{j}^{-1} + \sum_{\g\in\Gamma-\{e\}}
c_{\g}\g_{i}\,.$$
\item[{\emph{(R'2)}}] 
For any $u,v\in L$ and $i\neq j$:
$$[u_{i},v_{j}]= -\frac{k}{2} \sum_{\g\in\Gamma} \omega_{L}(\g u,v)
s_{ij}\g_{i}\g_{j}^{-1} \,.$$
\end{itemize}
\end{defi}
\epf

In this paper we will deal with the classification of  representations of the algebra $\h$ in which $x_i, y_i$  act by zero for any $i$, thus the relations we will actually use are the following:
\begin{itemize}
\item[{\emph{(R1)}}] 
For any $i\in [1,N]$: 
$$
0=  1+ \frac{k}{2} \sum_{j\neq i}\sum_{\g\in\Gamma}
s_{ij}\g_{i}\g_{j}^{-1} + \sum_{\g\in\Gamma-\{e\}}
c_{\g}\g_{i}\,.
$$
\item[{\emph{(R2)}}] 
For any $u,v\in L$ and $i\neq j$:
$$
0=\frac{k}{2} \sum_{\g\in\Gamma} \omega_{L}(\g u,v)
\g_{i}\g_{j}^{-1} \,,
$$
\end{itemize}
where, with abuse of notation, we wrote $\g_i$, $s_{ij}$ etc. \dots  for the images of the corresponding elements of $\h$ in a representation.

We end this section by recalling that in the rank $1$ case we have $\mathbf{\G}_1=\G$. In this case there's no action of $S_N$ and no parameter $k$, and the symplectic reflection algebra attached to $\G$ is simply:
$$
H_{1,c}(\G)=\frac{\C<x,y>}{I([x,y]=1+\sum_{\g\in\G-\{e\}}c_{\g}\g)} ,
$$
where $\;\C<x,y>=TL\;$, the tensor algebra of $\;L,\;$ and $\;I(-)\;$ stands for the ideal generated by the indicated relations. Representations of these algebras (deformations of Kleinian singularities) have been studied in \cite{CBH}.

It's easy to see that relations \emph{(R1)}, \emph{(R2)} in this case reduce to:
\begin{equation}\label{rank1}
0=1+\sum_{\g\in\G-\{e\}}c_{\g}\g.
\end{equation}

\subsection{Irreducible representations of wreath product groups}\label{irred}

For  reader's convenience and in order to introduce some important notation, we  recall the classification of irreducible representations for a wreath product group. Everything that follows is true for any finite group $\G$ and for representations over any algebraically closed field $F$ of characteristic $0$. For simplicity we will consider $F=\C$, the field of complex numbers. For complete proofs and details the reader should refer to \cite{JK}, Chapter $4$. 

A nice property of the wreath product group $\Gg$ is that the set of its irreducible representations $\mathrm{Irr}(\Gg)$  can be completely recovered from a knowledge of $\mathrm{Irr}(\G)$, using the representation theory of the symmetric group. 

Let $\{Y_1,\dots,Y_{\nu}\}$ denote the complete set of pairwise non-isomorphic representations of $\G$ over $\C$. Then a complete set of irreducible representations of $\G^N$ is given by $Y=Y_{h_1}\otimes\dots \otimes Y_{h_{N}}$ where  $(h_1,\dots,h_N)$ varies in $ [1,\nu]^N$.
If $N_h$ denotes the number of indices $i$ s.t.  $h_i=h$,  i.e. the number of factors of $Y$ equal to $Y_h$, $h=1,\dots,\nu$, then:
$$
\Nv=(N_1,\dots,N_{\nu})
$$ 
is called the \emph{type} of $Y$. 

We will say that two representations $Y$, $Y'$ are \emph{conjugate} if they have the same type.  This simply means that $Y=Y_{h_1}\otimes \dots \otimes Y_{h_N}$ and $Y'=Y_{h_{\s (1)}} \otimes \dots \otimes Y_{h_{\s(N)}}$ for some $\s\in S_N$, i.e. $Y'$ equals the representation $Y$ twisted by the outer automorphism of $\G^N$ that permutes the factors according to $\s$. It turns out that the role played by conjugate representations of $\G^N$ in recovering irreducible representations of $\Gg$ is exactly the same. This is essentially because, as one can easily argue from Definition \ref{wreath}, the outer automorphism induced by $\s\in S_N$ on $\G^N$ is a restriction of an inner automorphism in $\Gg$ (conjugation by the element $\s\in\Gg$).  So from now on we will consider only the representations of $\G^N$ that can be written as $Y=Y_1^{\otimes N_1}\otimes\dots\otimes Y_{\nu}^{\otimes N_{\nu}}$. Notice that the representations of this form are a complete set of irreducible, pairwise non-conjugate representations of $\G^N$.  

For any $h$,  we  denote by $S_{N_h}$ the subgroup of $S_N$ consisting of the permutations  that move only the indices $\{\sum_{i=1}^{h-1}N_i+1,\dots,\sum_{i=1}^h N_i\}$, corresponding to the  factors of $Y$ isomorphic to $Y_h$. We agree that  $S_{N_h}=\{1\}$ if $N_h=0$. Thus  we can consider the group:
$$
S_{\Nv}=S_{N_1}\times \dots\times S_{N_{\nu}}\subset S_N\subset S_N\ltimes \G^N $$
called the \emph{inertia factor} of $Y$.
Obviously any irreducible representation  $W$ of $S_{\Nv}$ is obtained as $W=W_1\otimes\dots\otimes W_{\nu}$, where $W_h$ is an irreducible representation of $S_{N_h}$.

The \emph{inertia subgroup} of  $Y$, instead, is defined to be:
$$
(\Gg)_Y= S_{\Nv}\ltimes \G^N\subset S_N\ltimes \G^N.
$$

Let's now consider  an irreducible representation of $\G^N$, $Y=Y_{1}^{\otimes N_1}\otimes\dots\otimes Y_{\nu}^{\otimes N_{\nu}}$.
There is a natural action of $(\Gg)_Y$  on $Y$ in which $\G^N$ acts  according to $Y$ and $S_{\Nv}$ permutes the factors. This representation can be shown to be irreducible. For simplicity we will keep the notation $Y$ for this representation.

Another easy way to obtain irreducible representations of  $(\Gg)_Y$ is extending an irreducible representation $W=W_{1}\otimes\dots\otimes W_{\nu}$ of $S_{\Nv}$ by making $\G^N$ act trivially. In this case  we will also  keep the notation $W$ for this extension.

Let's now consider the tensor product of $W$ and $Y$: 
$$
W\otimes Y=W_1\otimes\dots\otimes W_{\nu}\otimes Y_{1}^{\otimes N_1}\otimes\dots\otimes Y_{\nu}^{\otimes N_{\nu}}.
$$
Here  $S_{\Nv}$ acts both on $W$ and on $Y$( permuting the factors), while $\G^N$ acts only on $Y$. 
 This is  also an irreducible representation of $(\Gg)_Y$ (\cite{JK}, page 155). 

We can now obtain the induced representation of $\Gg$:
$$ 
W\otimes Y\uparrow:=\mathrm{Ind} _{\scriptscriptstyle{(\Gg)_Y}}^{\scriptscriptstyle{\Gg}}W\otimes Y 
$$
The following theorem holds.

\begin{thm}\label{wregroup}
The representation $W\otimes Y \uparrow$ is  irreducible and runs through a complete system of pairwise non-isomorphic irreducible representations of  $\Gg$ if $Y$ runs through a complete system of pairwise non-conjugate  irreducible representations of $\G^N$ and, while $Y$ remains fixed, $W$ runs through a complete system of pairwise non-isomorphic  irreducible representations of $S_{\Nv}$.
\end{thm}

In particular we have that, for a fixed $W$, the representation $W\otimes Y \uparrow$ depends only on the type of $Y$.   
With abuse of language we will call  \emph{type} of  $ W\otimes Y\uparrow$ the type of $Y$ as a representation of $\G^N$. We remark that the possible types of $W\otimes Y\uparrow$ are in bijection with the $\nu$-tuples $(N_1,\dots,N_{\nu})$, $N_h\geq 0$, $\sum_h N_h=N$ and that to  any such $\nu$-tuple  we can attach a proper partition of $N$, taking all the non-zero $N_h$s in $(N_1,\dots,N_{\nu})$ and ordering them in non-increasing order.

\rem\  A very easy  example is when  all the factors of $Y$ are the same , i.e. $Y=Y_h^{\otimes N}$ for some $h\in [1,{\nu}]$. The type of $Y$ is $(0,\dots 0,N,0,\dots 0)$ with an $N$ in the $h$th position and is associated to the partition of $N$ of  Young diagram a single row of length $N$, i.e. the partition corresponding to the trivial representation of $S_N$. For this reason we will call these representations of ``trivial type''.  In this case the inertia factor of $Y$ is $S_N$, its inertia subgroup coincides with $S_N\ltimes \G^N$, and we need no induction.  For any irreducible representation $W$ of $S_N$ we obtain the irreducible representation $W\otimes Y$ of  $S_N\ltimes\G^N$. 
 
\subsection{Irreducible representations for  $\G\subset SL(2,\C)$ and the McKay correspondence}\label{gammarepr}

In this section we will briefly recall some classical results about the finite subgroups of $SL(2,\C)$ and their representations that we will need in the next sections. 

It is well known (see for example \cite{Co}, Chapters $6$,$7$) that all the finite subgroups of $SL(2,\C)$, or equivalently the finite groups of quaternions, can be distinguished into two infinite series: 
\begin{itemize} 
 \item the cyclic groups $\mathcal{C}_{n+1}$ for any $n\geq 0$ ($\mathcal{C}_1=\{e\}$), of order $n+1$;
 \item the dicyclic groups $\mathcal{D}_{n-2}$ for $n\geq 4$, of order $4(n-2)$; 
 \end{itemize} 
  and three exceptional groups, that are the double coverings of the groups of rotations preserving regular polyhedra in ${\real}^3$ via the homomorphism of Lie groups $SU(2)\longrightarrow SO(3,\real)$:
  \begin{itemize}
 \item the binary tetrahedral group $\mathfrak{T}$, of order $24$;
  \item the binary octahedral group $\mathfrak{O}$, of order $48$;
 \item  the binary icosahedral group $\mathfrak{I}$, of order $120$.
  \end{itemize}
 
  The  terminology we used refers to the so called \emph{McKay correspondence}, as we are going to explain. In (\cite{MK})  McKay showed that \emph{``the eigenvectors of the Cartan matrices of affine type}  $\tilde{A}_n$, $\tilde{D}_n$, $\tilde{E}_6$, $\tilde{E}_7$, $\tilde{E}_8$ \emph{can be taken to be the columns of the character tables of the finite groups of quaternions''}. To  this end  he attached a graph  to any finite subgroup of $SL(2,\C)$  in the following way. Consider the set of irreducible non-isomorphic representations of a finite group $\G\subset SL(2,\C)$,  $I=\left\{Y_1,\dots ,Y_{\nu}\right\}$, and let $L$ be the tautological representation of $\G$, i.e. the natural representation of $\G$ as a subgroup of $SL(2,\C)$. Notice that $L$ is a self-dual representation.  Now build the graph  in which the set of vertices is $I$, and  the number of edges between two vertices $Y_h$ and $Y_{h'}$ is the multiplicity of the irreducible representation $Y_h$ in  $Y_{h'}\otimes L$ or equivalently, since $L$ is self-dual,  the multiplicity of  $Y_{h'}$ in  $Y_h\otimes L$. Any such graph turns out to be  an extended Dynkin graph with extending vertex corresponding to the trivial representation. If we label each vertex with the dimension of the corresponding representation the result is the following:

\vspace{0.8cm}  

\begin{figure}[h!]\label{cyclicfig}
\center{\includegraphics{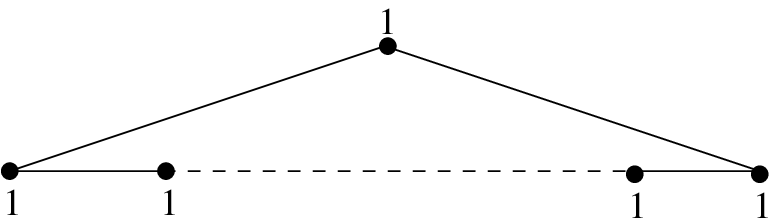}}
\caption{$\tilde{A}_n$ for $\G=\mathcal{C}_{n+1},\: n\geq 0$ }
\end{figure}

\vspace{0.8cm}

\begin{figure}[h!]\label{dicyclicfig}
\center{\includegraphics{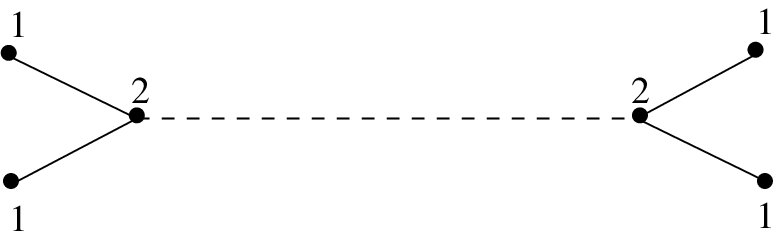}}
\caption{$\tilde{D}_n$ for $\G=\mathcal{D}_{n-2},\: n\geq 4$}
\end{figure}

\vspace{0.8cm}

\begin{figure}[h!]
\center{\includegraphics{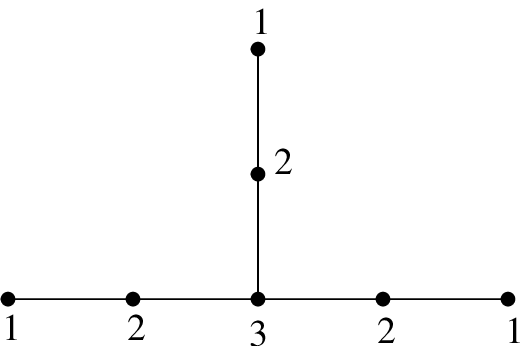}}
\caption{$\tilde{E}_6$ for $\G=\mathfrak{T}$ }
\end{figure}

\vspace{0.8cm}

\begin{figure}[h!]
\center{\includegraphics{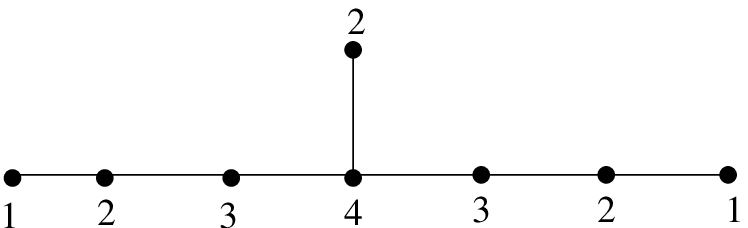}}
\caption{$\tilde{E}_7$ for $\G=\mathfrak{O}$ }
\end{figure}

\vspace{0.8cm}

\begin{figure}[h!]
\center{\includegraphics{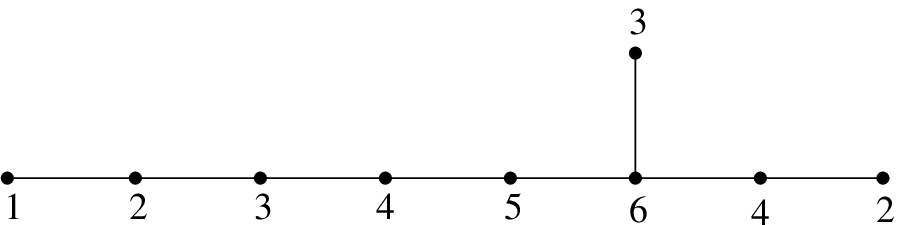}}
\caption{$\tilde{E}_8$ for $\G=\mathfrak{I}$ }
\end{figure}

\vspace{0.8cm}

Notice that the adjacent vertices to a fixed vertex $Y_h$ correspond to the irreducible components of the representation $Y_h\otimes L$.

\subsection{Representations of $S_N$ with rectangular Young diagram}\label{rectangle}

 In what follows we will use the following standard results from  representation theory of the symmetric group that we will state without proof. Denote by $\mathfrak{h}$ the reflection representation of $S_N$. For a Young diagram $\lambda$ we denote by $W_{\lambda}$ the corresponding irreducible
representation of $S_N$ and by $C(\lambda)$ the content of $\lambda$, i.e. the sum of signed distances of the cells from the diagonal. 

\begin{lem}\label{rectangular}
\begin{center}
\begin{itemize}
\item[{i)}] $\Hom_{S_N}(\mathfrak{h}\otimes W_{\lambda},W_{\lambda})=\C^{m-1}$, where $m$ is the number of corners of the Young diagram $\lambda$. In particular  $\Hom_{S_N}(\mathfrak{h}\otimes W_{\lambda},W_{\lambda})=0$ if and only if $\lambda$ is a rectangle;
\item[{ii)}] the element $C=s_{12}+s_{13}+ \cdots +s_{1N}$ acts by a scalar in $W_{\lambda}$ if and only if $\lambda$ is a rectangle. In this case $C|_{W_{\lambda}}=\frac{2\,C(\lambda)}{N}$;
\item[{iii)}] if $\lambda$ is a rectangular Young diagram of height $a$ and width $b$, then: $C(\lambda)=\frac{(b-a)\,N}{2}$.
\end{itemize}
\end{center}
\end{lem}

\epf

\section{The main theorem}\label{main}

 Our main theorem classifies the irreducible representations of $\h$ for values of $(k,c)$ with $k\neq 0$, in which the generators $x_i$, $y_i$ act by $0$ for any $i=1,\dots,N$. Since such representations can be considered as irreducible representations of $\Gg$, the problem reduces to classifying the irreducible representations of $\Gg$ that can be  extended to a representation of the associated symplectic reflection algebra with trivial action of the generators $x_i$, $y_i$  for such values of $(k,c)$. On the other hand, we will  show at the beginning of the next section that if an irreducible representation of $\Gg$ extends to a representation of $\h$  then $x_i$, $y_i$ must necessarily act by $0$. 

For $\G=\{e\}$ it is easy to see that the  algebra $H_{1,k,c}(S_N)$ has no finite dimensional representations. In fact  $H_{1,k,c}(S_N)$   always contains a copy of the Weyl algebra (generated by the elements $x_1+\dots+x_N$, $y_1+\dots+y_N$)  that has no finite dimensional representations. We will thus consider the case $\G\neq\{e\}$.
Before stating the theorem we need to introduce some notation:
\begin{itemize}
\item $\nu$ will denote the number of conjugacy classes $\{C_1,\dots,C_{\nu}\}$  of $\G$ , with $C_1=\{e\}$ , $|C_s|$ will be  the cardinality of the class $C_s$, and $c_s$  the value of the class function $c$ on $C_s$;
\item for any irreducible representation $Y_h$ of $\G$, $\chi_{Y_h}(C_s)$ will be the value of the character of $Y_h$  on the class   $C_s$ . 
\end{itemize}

With this notation, the complex  number $ \frac{|C_s|\,\chi_{Y_h}(C_s)}{\di Y_h}$ is  the scalar corresponding to the central element  $\sum_{\g\in C_s}\g$ in the irreducible representation $Y_h$.

\begin{thm}\label{mainth}
Let $\G\neq\{e\}$. 
\begin{enumerate}
\item[I)] For $k\neq 0$ an  irreducible  representation $W\otimes Y\uparrow$ of $\Gg$ of type $(N_1,\dots, N_{\nu})$ extends to a representation of the associated symplectic reflection algebra $\h$ if and only if the following two conditions are satisfied:

\begin{enumerate}
\item[i)] $W=W_1 \otimes\dots\otimes W_{\nu}$, where $W_h$ is an irreducible representation of $S_{N_h}$ with rectangular Young diagram $\lambda_h$, of size $a_h\times b_h$;

\item[ii)] for any $h\neq h'$ s.t. $N_h,N_{h'}\neq 0$, $\Hom_{\G}(Y_h\otimes L ,Y_{h'})=0$, where $L$ is the tautological representation of $\G$. In other words,  any two  non-isomorphic representations $Y_h$, $Y_{h'}$ of $\G$ occurring in the type of $Y$  must  be  non-adjacent vertices  in the extended  Dynkin diagram attached to $\G$.
\end{enumerate}

\item[II)] $W\otimes Y\uparrow$ extends  for any value of the parameter $(k,c)$ belonging to the intersection of the hyperplanes:
\begin{equation}\label{maintwo}
\mathcal{H}_h:\qquad \di Y_h+(b_h-a_h)\,\frac{k}{2}\,|\G|+ \sum_{s=2}^{\nu}c_s\,|C_s|\,\chi_{Y_h}(C_s)=0 
\end{equation} 
for any $h\in\{1,\dots,\nu\}$ s.t. $N_h\neq 0$,  i.e. for any representation $Y_h$ occurring in the type of $Y$.
The space of the solutions of this system of equations has dimension $\nu-r$ where $r=\#\{N_h \mbox{\ s.t.\ } N_h\neq 0\}$. 
\end{enumerate}
\end{thm}

The next section is  devoted to the proof of Theorem \ref{mainth}.

\section{Proof of Theorem \ref{mainth}}

As stated at the beginning of the previous section we start by proving the following:
\begin{thm}\label{zero}
Let $M$ be a $H_{1,k,c}(\Gg)$-module such that $M|_{\Gg}$ is irreducible. Then the generators $x_i$, $y_i$ act by zero on $M$ for any $i=1,\dots,N$.
\end{thm}

\pf\, Without loss of generality  consider the elements $x_1,\,y_1 \in H_{1,k,c}(\Gg)$. From Section \ref{symal} we know these elements commute with the elements $\g_i$ for $i\neq 1$, and the action of $\g_1$ by conjugation on such elements  corresponds to the action of $\g$ on the basis vectors $x$, $y$ respectively in the tautological representation $L$ of $\Gamma$. Thus we can view $x_1$, $y_1$ as a basis for the representation:
$$
L\otimes\underbrace{\C\otimes\dots\otimes\C}_{N-1}
$$
of $\Gamma^N$, where $\C$ is the trivial one-dimensional representation.  So we have that the action of $x_1$, $y_1$ on $M$ induces maps of $\Gamma^N$-modules:
$$
\left(L\otimes\C\otimes\dots\otimes\C\right)\otimes M\longrightarrow M.
$$  

 But now from Section \ref{irred}  we have that, as a $\Gamma^N$-module, $M$ decomposes in irreducibles as:
$$
\bigoplus_{\s}Y_{h_{\s(1)}}\otimes\dots\otimes Y_{h_{\s(N)}}
$$
where $\s$ are permutations in $S_N$ and factors  may appear with some multiplicity. Thus composing with the $\Gamma^N$-module maps given by the injections and projections of the direct factors we have that $x_1$, $y_1$ induce $\Gamma^N$-module maps:
$$
\left(L\otimes Y_{h_{\s(1)}}\right)\otimes\dots\otimes Y_{h_{\s(N)}}\longrightarrow Y_{h_{\s'(1)}}\otimes\dots \otimes Y_{h_{\s'(N)}}
$$
for any $\s$, $\s'$. Since the $Y_h$s are irreducible $\Gamma$-modules, in order for such a map to be non-zero, we  must have  $ Y_{h_{\s(i)}}\cong  Y_{h_{\s'(i)}}$ for any $i\geq 2$. This implies $ Y_{h_{\s(1)}}\cong Y_{h_{\s'(1)}}$ and we get a homomorphism:
$$
 L\otimes Y_{h_{\s(i)}} \longrightarrow Y_{h_{\s(i)}}.
$$
But  such a homomorphism must be zero as explained in Section \ref{gammarepr},  as  extended Dynkin diagrams have no loop-vertices. 
We deduce $x_1$, $y_1$ act trivially on $M$.
\vspace{0.2 cm}

\hspace{14.2 cm} $\Box$
\vspace{0.1 cm} 

To prove Theorem \ref{mainth} we will use simple classical results  from the representation theory of finite groups. From now on we will assume $\G\neq\{e\}$.

\subsection{The relations \emph{(R2)}}\label{rela2}

It turns out that the relations \emph{(R2)} have an easy interpretation in terms of the extended  Dynkin diagram attached to the group $\G$ in the McKay correspondence. Let $L$ be the tautological representation of $\G$. We have the following theorem.

\begin{prop}\label{r2}
If $W \otimes Y\uparrow$ is a representation of $\Gg$ of type $(N_1,\dots,N_{\nu})$, then the operators of the corresponding matrix representation  satisfy \emph{(R2)} for $k\neq 0$ if and only if for any pair $h,h'$ s.t. $N_h, N_{h'}\neq 0$, $\Hom_{\G}(L\otimes Y_{h},Y_{h'})=0$ i.e. if and only if  $Y_h$, $Y_{h'}$  are not adjacent vertices of the extended Dynkin diagram associated to $\G$ in the McKay correspondence. 
\end{prop} 
\pf\  Relations \emph{(R2)} are satisfied for $k\neq 0$ if and only if 
$$
\sum_{\g\in\G}\omega_L(\g u,v)\,\g_i\g_j^{-1}=0\qquad \forall\, u,v\in L,\quad \forall\, i\neq j\, .
$$ 
We observe that, for any $Y$, the subgroup $\G^N$ is contained in the inertia subgroup of $Y$ and is normal in $\Gg$. For this reason the induced representation $W\otimes Y\uparrow$ can be written as:
\begin{equation}\label {direct}
\s_1\cdot(W\otimes Y)\oplus\dots\oplus \s_{\M}\cdot(W\otimes Y)
\end{equation}
where  $M=\frac{N!}{N_1!\dots N_{\nu}!}$, and $\{\s_1,\dots,\s_{\M}\}$  is a set of representatives of the  left cosets of the inertia factor $(\Gg)_Y$ in $\Gg$, that can be chosen to be all in $S_N$. The action of an element $g\in \Gg$ on a vector $\s_l\cdot v$ is defined as follows:
$$
 g(\s_l\cdot v)= \s_r\cdot(g'v)\qquad \mbox{where\ }  g\s_l=\s_rg'\qquad g'\in (\Gg)_D\, .
$$
By the normality of  $\G^N$, all  the direct factors of (\ref{direct}) are stable under the action of  $\sum_{\g\in\G}\omega_L(\g u,v)\,\g_i\g_j^{-1}$, thus  this operator  has a block diagonal form. The $l$th block   corresponds to the operator   $A(s,t)=\sum_{\g\in\G}\omega_L(\g u,v)\,\g_s\g_t^{-1}$, with $(s,t)=(\s_l^{-1}(i),\s_l^{-1}(j))$,
 in the representation $W\otimes Y$ of $(\Gg)_Y$. We are reduced now to show that any such block is zero if and only if the conditions of the proposition  are satisfied. Since  the action of  $A(s,t)$ is trivial on $W$, we can suppose $W$ to be trivial $1$-dimensional, thus  $W\otimes Y\cong Y$. Without loss of generality we can also suppose $s\leq t$. Since the bilinear form $\omega_L$ is non degenerate and $u$, $v$ vary in all $L$ we have:
$$
A(s,t)=0 \Leftrightarrow   \sum_{\g\in G}\g\otimes\g\otimes\g^{-1}|_{L\otimes Y_{h_s}\otimes Y_{h_t}}=0\Leftrightarrow  \sum_{\g\in G}\g\otimes\g\otimes {\g^{-1}}^{\ast}|_{L\otimes Y_{h_s}\otimes Y_{h_t}}=0
$$
where ``$\,^{\ast\,}$'' denotes the transposition. Now, if we denote by $Y_{h_t}^{\ast}$ the dual representation of $Y_{h_t}$, we notice that the last operator corresponds to the operator:
$$
 \sum_{\g\in G}\g\otimes\g\otimes\g|_{L\otimes Y_{h_s}\otimes Y^{\ast}_{h_t}} , 
$$
that is a multiple of the projector on the invariants of the representation $L\otimes Y_{h_s}\otimes Y_{h_t}^{\ast}$.

Now this is zero if and only if $\Hom_{\G}(L\otimes Y_{h_t}, Y_{h_s})=0$.  
From Section \ref{gammarepr} we know  this happens exactly when $Y_{h_s}$, $Y_{h_t}$ are non adjacent vertices in the Dynkin diagram attached to $\G$.   
\epf

\rem\ Notice that  when $Y$ is of trivial type $(0,\dots,N,\dots 0)$,  i.e. when all the factors of $Y$ are the same, Proposition \ref{r2} implies that the conditions \emph{(R2)} are automatically  satisfied (Dynkin diagrams corresponding to non-trivial finite subgroups of $SL(2,\C)$ have no loop vertices).  

\subsection {The relations \emph{(R1)}}\label{rela1}
The only thing we are left to do now is analyzing the conditions for relations \emph{(R1)} to be satisfied. We will begin from the easiest case of $W\otimes Y\uparrow$ with $Y$ of type $(0,\dots,N,\dots,0)$, corresponding to the trivial partition of $N$ (a single row with $N$ boxes). We have the following proposition.

\begin{prop}\label{cycsquare}
For  $k\neq 0$ a representation $W\otimes Y$ of $\Gg$ of trivial type $(0,\dots,N,\dots  0)$  extends to a representation of $\h$ if and only if  the following conditions are satisfied:
\begin{enumerate}
\item[i)] the representation $W$ of $S_N$ corresponds to a rectangular Young diagram;
\item[ii)] the parameter  $(k,c)$   satisfies the corresponding equation in part $II)$ of Theorem \ref{mainth}.
\end{enumerate}
\end{prop}
\pf\ As we observed in the previous subsection, Proposition \ref{r2} implies that in this case relations \emph{(R2)} are satisfied. Thus we only have to consider relations \emph{(R1)}. 

We will begin with an easy clarifying example. When  $Y=Y_h^{\otimes N}$ with $Y_h$ one dimensional, it is straightforward  to check that Proposition \ref{cycsquare} holds. In this case in fact, the permutation action of $S_N$ on $Y_h^{\otimes N}$ is trivial and $W\otimes Y_h^{\otimes N}\cong W$ as $S_N$-modules. Thus the  relation \emph{(R1)} for a fixed $i$ looks like:
$$
0=\frac{k}{2}\,|\G|\, \sum_{j\neq i}s_{ij}=-1-\sum_{\g\in\Gamma-\{e\}}
c_{\g}\,\chi_{Y_h}(\g)= -1-\sum_{s=2}^{\nu}c_{s}\,|C_s|\,\chi_{Y_h}(C_s)\, ,
$$
where $\chi_{Y_h}(C_s)$ is the value of the character of $Y_h$ on the conjugacy class $C_s$. For $k\neq 0$ we  have:
 \begin{equation}\label{onedi}
 \sum_{j\neq i}s_{ij}=\frac{  2\left(-1- \sum_{s=2}^{\nu}c_{s}\,|C_s|\,\chi_{Y_h}(C_s)\right)}{k\,|\G|}\, .
\end{equation}
So $ \sum_{j\neq i}s_{ij}$ must act as a scalar and Lemma (\ref{rectangular}) part $ii)$ implies that $W$ must have rectangular Young diagram $\lambda$ of some size $a\times b$. We remark that Lemma \ref{rectangular} part $ii)$  $iii)$  implies that the element $\sum_{j\neq i}s_{ij}$  acts as the scalar   $\frac{2\,C(\lambda)}{N}=(b-a)$ in this representation. Substituting this value in equation (\ref{onedi}) we get the result in the $1$-dimensional case.  Notice that this first consideration solves completely   the case when $\G$ is cyclic.

Let's now suppose $\di\,W=m$ and $\di\,Y_h=n>1$. 
 We rewrite relations \emph{(R1)} as follows:
\begin{equation}\label{constant}
 - 1- \sum_{\g\in\Gamma-\{e\}} c_{\g}\g_{i}= \frac{k}{2} \sum_{j\neq i}\sum_{\g\in\Gamma}s_{ij} \g_{i}\g_{j}^{-1}\,.
\end{equation}
 We observe that the LHS of (\ref{constant}) is a central element of the group algebra $\C[\G^N]$  (due to the fact  that $c$ is a class function), and that, as $\G^N$-module, $W\otimes Y$  is isomorphic to a direct sum of $\,\di W\,$ copies of the irreducible representation $Y_h^{\otimes N}$. Thus the LHS acts as a scalar in this representation. More precisely we have:
$$
\sum_{\g\in\Gamma-\{e\}} c_{\g}\g_{i}=\sum_{s=2}^{\nu}\frac{c(s)\,|C_s|\,\chi_{Y_{h}}(C_s)}{\di Y_h}\,.
$$ 
So we must have that   $ \frac{k}{2} \sum_{j\neq i}\sum_{\g\in\Gamma} s_{ij}\g_{i}\g_{j}^{-1}$ is a scalar. We will show that this operator has a block form that reduces equation (\ref{constant}) to equation (\ref{onedi}).

For this let's  take  any two bases $\{v_1,\dots,v_n\}$, $\{w_1, \dots,w_{m}\}$ for  $Y_h$ and  $W$ respectively. The vectors $\{ \V_{I}=v_{i_1}\otimes\dots\otimes v_{i_N}\}$, where the multi-index  $I=(i_1,\dots,i_N)$ varies in $[1,n]^{N}$, are clearly a basis  for $Y=Y_h^{\otimes N}$. We can give the multi-indices $I$ a total ordering $I_1,\dots,I_{n^N}$ using  the lexicographic order.  Consider now the basis of $W\otimes Y_h^{\otimes N}$ given by the vectors: 
$$
Z_1=w_1\otimes \V_{I_1},\,\dots\,, Z_m=w_m\otimes \V_{I_1},\,\dots\, ,  Z_{n^N}=w_1\otimes \V_{I_{n^{\N}}},\,\dots \, , Z_{{m n^N}}=w_{m}\otimes \V_{I_{n^{\N}}}\, .
$$
Any transposition  $s_{ij}\in S_N$ induces a permutation $\tilde{s}_{ij}$ (of order $2$) on the set $\{I_1,\dots,I_{n^N}\}$ thus on the vectors of the basis $\{Z_1,\dots,Z_{mn^N}\}$. Let's now  denote by $A_W(s_{ij})$ the operator (of size $m\times m$) for $s_{ij}$ in the representation $W$, and by $O_m$ the $0$-operator of size $m\times m$. It is easy to see that, using the basis $\{Z_1,\dots,Z_{mn^N}\}$, we can obtain a block form for the operator  $s_{ij}$ in the representation $W\otimes Y$ from the block diagonal operator
$$
\left(\begin {array}{cccc}
A_{W}(s_{ij})&O_m &\cdots  & O_m\\
O_m &\ddots &\cdots & O_m\\
\vdots &\vdots & \ddots& \vdots\\
O_m &\cdots &O_m & A_{W}(s_{ij})
\end{array}\right)
$$
 by simply permuting the columns according to $\tilde{s}_{ij}$.
Using this, we can compute a block form for $s_{ij}\sum_{\g\in\Gamma}\g_i\g_j^{-1}$. We denote each block, of size $m\times m$, by its position  $(\R,\T)$, where $\R=(r_1,\dots,r_N)$, $\T=(t_1,\dots,t_N)$ are multi-indices. 
We have the following formulas for the blocks :
\begin{itemize}
\item for $(\R,\T)$ with $\R$ differing from $\T$ at most for the pair  of indices $(r_i, r_j)$:
$$
(\R,\T)=A_W(s_{ij})\,\sum_{\g\in G}\alpha_{r_j t_i}(\g)\,\alpha_{r_i t_j}({\g}^{-1})\, ,
$$
where $\alpha_{r_j,t_i}(\g)$ are the matrix coefficients in the representation $Y_h$;\\
\item
for $(\R,\T)$ with $\R$ differing from $\T$ for indices different from $r_i, r_j$
$$
(\R,\T)=O_m\, .
$$
\end{itemize}
Summing up over $j\neq i$ we can now rewrite relations \emph{(R1)} in block form  for each $i\in[1,N]$:
\begin {enumerate}
\item  for $(\R,\T)$ with  $\R$ differing from $\T$  at most for the index $r_i$ 
 \begin{equation}\label{uno}
\frac{k}{2}\,\sum_{j\neq i}A_{W}(s_{ij})\,\sum_{\g\in\G} \alpha_{r_j t_i}(\g)\alpha_{r_i r_j}({\g}^{-1})=-\delta_{r_i t_i}\,I_m\,\left(1+\sum_{s=2}^{\nu}\frac{c_{s}\, |C_s|\,\chi_{Y_h}(C_s)}{\di Y_h}\right)
\end{equation}
where 
$\delta_{r_it_i}=\left\{\begin{array}{ll}
0 & \textrm{if $r_i\neq t_i$}\\ 1 & \textrm{if $r_i=t_i$}
\end{array}\right.$;\\
\item  for $(\R,\T)$ with  $\R$ differing from $\T$  at  least for an  index $r_j$, $j\neq i$, and at most for the pair  of indices $(r_i,r_j)$ 
\begin{equation}\label{due}
\frac{k}{2}\,A_{W}(s_{ij})\,\sum_{\g\in\G} \alpha_{r_j t_i}(\g)\alpha_{r_i t_j}({\g}^{-1}) =0\, .
\end{equation}
\end{enumerate} 
In all the other cases we only obtain  trivial relations.

Now we  observe   that, using the orthogonality property of matrix coefficients of  irreducible representations of a finite group, we get:
 $$
\sum_{\g\in \G}\alpha_{r_j,t_i}(\g)\,\alpha_{r_i,t_j}({\g}^{-1})=\delta_{r_it_i}\,\delta_{r_jt_j}\,\frac{|\G|}{\di Y_h}\,. 
$$
Substituting these values in  equation (\ref{due}) we obtain  trivial relations.  From equation (\ref{uno}), instead, we obtain that $\sum_{j\neq i}A_W(s_{ij})$ must be a scalar operator. Thus Lemma \ref{rectangular} implies that the Young diagram $\lambda$ attached to $W$ is a rectangle, of some size $a\times b$, and that $\sum_{j\neq i}A_W(s_{ij})$ acts on $W$ as the scalar $(b-a)$.  
Thus from  equation (\ref{uno})  we obtain the equation:
$$
\di Y_h+(b-a)\,\frac{k}{2}\,|\G|+ \sum_{s=2}^{\nu}c_s\,|C_s|\,\chi_{Y_h}(C_s)=0 
$$
which is exactly the  equation (\ref{maintwo}) for the hyperplane $\mathcal{H}_h$ in  Theorem \ref{mainth}, part $II)$. Notice that, in this case, we get a single equation since $Y_h$ is the only factor appearing in $Y$.
\epf
\bigskip

We will now  analyze the   cases when  the inertia factor of $Y$ is not the entire $S_N$ and an actual induction is needed to build the representation $W\otimes Y\uparrow $. If the type of $Y$ is $\Nv=(N_1,\dots,N_{\nu})$, then  the inertia factor is $S_{\Nv}=S_{N_1}\times\dots\times S_{N_{\nu}}$ and we have:
$$
W\otimes Y\uparrow = \s_1\cdot\left( W\otimes Y\right)\oplus\dots\oplus \s _{\M}\cdot \left(W\otimes Y\right)
$$
where $M=\frac{N!}{N_1!\cdots N_{\nu}!}$ and $\{\s_1,\dots,\s_{\M}\}$ is a set of representatives for the left cosets of $S_{\Nv}$ in $S_N$.
\vspace{0.3 cm}

\rem \  Let's denote by $[\s]$ the left coset of $\s$ with respect to $S_{\Nv}$. An easy computation shows that for any transposition $s_{ij}$ and any permutation $\s$: 
$$ 
[s_{ij}\s]=[\s] \Leftrightarrow s_{\s^{-1}(i)\s^{-1}(j)}\in S_{\Nv}\, .
$$
Moreover we observe that for any $\s\in S_N$ and any $i=1,\dots,N$:
$$
\g_i\s=\s\g_{\s^{-1}(i)}\, .
$$
We are now ready to prove the following result.

\begin{prop}\label{nontrivialtype}
 For  a representation $W\otimes Y\uparrow$ of $\Gg$ of non-trivial  type $(N_1,\dots ,N_{\nu})$ relations \emph{(R1)} are satisfied for some non-zero values of $k$ if and only if:
\begin{enumerate} 
\item[i)] $W=W_1\otimes\dots \otimes  W_{\nu}$, with $W_h$ irreducible representation of $S_{N_h}$ with rectangular Young diagram; 
\item[ii)] the parameter $(k,c)$   satisfies the corresponding system of equations in part $II)$ of Theorem \ref{mainth}. 
\end{enumerate}
\end{prop}

\pf\ Let $Y$ be a representation of $\G^N$ of  type $(N_1,\dots,N_{\nu})$. 
We observe that for any $W$, if we choose $\{\s_1,\dots,\s_M\}\subset S_N$  representatives of the left cosets of $(\Gg)_Y$ in $\Gg$:
\begin{equation}\label{decomposition}
W\otimes Y\uparrow = \s_1\cdot\left( W\otimes Y\right)\oplus\dots\oplus \s _{\M}\cdot \left(W\otimes Y\right),
\end{equation}
 is a $\G^N$-stable decomposition of $W\otimes Y\uparrow$.  For any representative $\s_l$, let's  denote by $\s_l\, Y$  the representation of $\G^N$ with same underlying vector space as $Y$ and with the action on $Y$ twisted by the automorphism induced by ${\s_l}$ on $\G^N$  (the action of $\g_i$ on $\s_l\,Y$ is the same as the action of $\g_{\s^{-1}(i)}$ on $Y$). Since $\G^N$ acts trivially on $W$, as a $\G^N$-module the subspace  $\s_l\cdot\left( W\otimes Y\right)$ is isomorphic to a direct sum of copies of the irreducible representation  $\s_l\, Y$.  So, for a fixed $i$, the $\G^N$-central operator $- 1- \sum_{\g\in\Gamma-\{e\}} c_{\g}\g_{i}$ preserves the subspaces $\s_l\cdot\left( W\otimes Y\right)$ and acts as a scalar on each of them. For any vector $ \s_l\cdot v\in\s_l(W\otimes Y)$ we have:
\begin {equation}\label{action1}  
\left(- 1- \sum_{\g\in\Gamma- \{e\}} c_{\g}\g_{i}\right)(\s_l\cdot v)=\s_l\cdot (- v)+\s_l\cdot \left(\left(- \sum_{\g\in\Gamma-\{e\}} c_{\g}\g_{\s_l^{-1}(i)}\right)v\right)= C (\s_l\cdot v).
\end{equation}

The action of $\frac{k}{2} \sum_{j\neq i}\sum_{\g\in\Gamma} s_{ij}\g_{i}\g_{j}^{-1}$  on such a vector  is instead:
\begin{equation}\label{action2}
\left(\frac{k}{2} \sum_{j\neq i}\sum_{\g\in\Gamma} s_{ij}\g_{i}\g_{j}^{-1}\right)(\s_l\cdot v)=\sum_{j\neq i}\s_{r(ijl)}\cdot \left(\left(\frac{k}{2}\sum_{\g\in\Gamma}\tilde{\s}_{ijl}\g_{\s_l^{-1}(i)}\g_{\s_l^{-1}(j)}^{-1}\right)v\right)
\end{equation}
where $\s_{r(ijl)}$ is the representative in the set $\{\s_1,\dots,\s_M\}$ of the coset $[s_{ij}\s_l]$ and $\tilde{\s}_{ijl}\in S_{\Nv}$ is the unique element s.t. $s_{ij}\s_l=\s_{r(ijl)} \tilde{\s}_{ijl}$. 

Relations \emph{(R1)} are satisfied if and only if these two actions are the same. In particular $\frac{k}{2} \sum_{j\neq i}\sum_{\g\in\Gamma} s_{ij}\g_{i}\g_{j}^{-1}$ must preserve the subspace $\s_l\cdot(W\otimes Y)$.
But let's  now look at equation (\ref{action2}) and take  $j\neq i$ s.t.   $[s_{ij}\s_l]=[\s_r]$, $r\neq l$ i.e. $s_{ij}$  ``moves''  the subspace  $\s_l\cdot(W\otimes Y)$ sending it to  the subspace $\s_r\cdot(W\otimes Y)$.  Then we  have: $s_{\s_l^{-1}(i) \s_l^{-1}(j)}\notin S_{\Nv}$. This means that  the representations $Y_{\s_l^{-1}(i)}$,  $Y_{\s_l^{-1}(j)}$ are not isomorphic. As a consequence, arguing as in Section \ref{rela2}, we have that 
$$
\sum_{\g\in\G}\g_{\s_l^{-1}(i)}\g_{\s_l^{-1}(i)}^{-1}=0
$$ 
in the representation $W\otimes Y$, hence $s_{ij}$ sends the subspace $\s_l\cdot(W\otimes Y)$ to $0$. This means that $\frac{k}{2} \sum_{j\neq i}\sum_{\g\in\Gamma} s_{ij}\g_{i}\g_{j}^{-1}$ indeed preserves the subspace $\s_l\cdot(W\otimes Y)$ and that relations \emph{(R1)} split up into equations that can be checked on the subspaces $\s_l\cdot(W\otimes Y)$ . So in equation (\ref{action1}) it is enough to take the sum over the $j\,$s s.t. $[s_{ij}\s_l]=[\s_l]$.  
Moreover we know that if $[s_{ij}\s_l]=[\s_l]$ then $s_{ij}\s_l=\s_l\, s_{\s_l^{-1}(i) \s_l^{-1}(j)}$ i.e. $\tilde{\s}_{ijl}= s_{\s_l^{-1}(i) \s_l^{-1}(j)}$. 
Hence, for a fixed $i$, if $\s_l^{-1}(i)=p$ the relations \emph{(R1)} reduce to the following equations:

\begin{equation}\label{finaleq}
 \frac{k}{2} \sum_{\begin{array}{c} q\neq p\\
                                  s_{p q}\in S_{N_{h_p}}
\end{array}}
s_{p q }\sum_{\g\in\Gamma}\g_{p}\g_q^{-1}=-1-\sum_{\g\in\G-\{e\}}c_{\g}\g_{p}
\end{equation}
where the identity must be considered in the representation $W\otimes Y$ of $S_{\Nv}\ltimes \G^N$ and $p\in\{\s_1^{-1}(i),\dots,\s_{\M}^{-1}(i)\}$. For any  $p$, equation (\ref{finaleq}) is exactly the $p$th equation   of relations \emph{(R1)} for  the extension of the representation of trivial type $W_{h_p}\otimes Y_{h_p}^{\otimes N_{h_p}}$ of $S_{N_{h_p}}\ltimes \G^{N_{h_p}}$ to the algebra $H_{1,k,c}(S_{N_{h_p}}\ltimes \G^{N_{h_p}})$. It is easily checked that, letting $i$ and $\s_l$ vary,  we  obtain all the relations for the extension of the representations $W_h\otimes Y_h^{\otimes N_h}$ of $S_{N_h}\otimes \G^{N_h}$ for any  $N_h\neq 0$ . Using Proposition \ref{cycsquare} we get the result.

\epf

\subsection{The conditions on the parameter $(k,c)$}\label{parameter}

Now that we found out which representations of $\Gg$ can potentially be extended to representations of $\h$ for  $k\neq 0$ we would  like to show that such extensions exist for a non-empty set of values of $(k,c)$.  This  amounts to prove that the system of equations in Theorem \ref{mainth} part $II)$ admits solutions. Fix a representation $W\otimes Y\uparrow$ of $\Gg$ of type $(N_1,\dots,N_{\nu})$ satisfying conditions $i)$ $ii)$ of Theorem \ref{mainth} part $I)$. We have the following proposition.

\begin{prop}\label{system}
 If  $r=\#\{N_h\mbox{\ s.t.\ }N_h\neq 0\}$, then the space of the solutions for the system of equations in part $II)$ of Theorem \ref{mainth} has  dimension $\nu-r$. 
\end{prop}

\pf\ 
By condition $ii)$ we have that: 
$$
r=\# \{N_h \mbox{\ s.t.\ } N_h\neq 0\}< \nu=\#\{\mbox {vertices in the extended Dynkin diagram of \ } \G\}\,.
$$
Without loss of generality we can suppose $N_1,\dots ,N_r\neq 0 $, $N_h=0$ for $h>r$. So in  matrix form the system has size $r\times\nu$, with $r\leq\nu-1$:
$$
\left(\begin{array}{cccc}\frac{(b_1-a_1)\,|\G|}{2} & |C_2|\chi_{Y_1}(C_2) &\dots &|C_{\nu}|\chi_{Y_1}(C_{\nu})\\
\frac{(b_2-a_2)\,|\G|}{2} & |C_2|\chi_{Y_2}(C_2) &\dots &|C_{\nu}|\chi_{Y_2}(C_{\nu})\\
\vdots &\vdots &\vdots &\vdots\\
 \frac{(b_r-a_r)\,|\G|}{2} & |C_2|\chi_{Y_r}(C_2) &\dots &|C_{\nu}|\chi_{Y_r}(C_{\nu})\end{array}\right) 
\left(\begin{array}{c} k\\c_2\\ \vdots\\c_{\nu}\end{array}\right)=\left(\begin{array}{c} -dim  Y_1\\ -dim Y_2\\ \vdots\\ -dim Y_r\end{array}\right)\, . 
$$
But now we have: 
$$
\rk \left( \begin{array}{ccc} |C_2|\chi_{Y_1}(C_2) &\dots &|C_{\nu}|\chi_{Y_1}(C_{\nu})\\
 |C_2|\chi_{Y_2}(C_2) &\dots &|C_{\nu}|\chi_{Y_2}(C_{\nu})\\
\vdots &\vdots &\vdots\\
|C_2|\chi_{y_r}(C_2) &\dots &|C_{\nu}|\chi_{Y_r}(C_{\nu})\end{array}\right)=\rk  \left( \begin{array}{ccc} \chi_{Y_1}(C_2) &\dots &\chi_{Y_1}(C_{\nu})\\
 \chi_{y_2}(C_2) &\dots &\chi_{Y_2}(C_{\nu})\\
\vdots &\vdots &\vdots\\
\chi_{Y_r}(C_2) &\dots &\chi_{Y_r}(C_{\nu})\end{array}\right)=r\, .
$$

In fact the rows $R_1,\dots ,R_r$ on the RHS are rows of the character table for $\G$ from which we have erased the entries $\chi_{Y_h}(e)=\di Y_h$. If a non-trivial linear combination $\sum_{h=1}^ra_h R_h$ of these rows is zero then  the class  function $X=\sum_{h=1}^{r}a_h\chi_{Y_h}$   satisfies the equation:  
 $$
X(\g)=0,\quad \forall\,\g\in\G-\{e\}\, .
$$
 This is possible  only  if $X=m \rho$, where $m\in\C$ and  $\rho$ is the character of the regular representation. Now we must have $m\neq 0$ since characters of non-isomorphic irreducible representations are linearly independent.  But $m\neq 0$  is also  impossible since, by condition $ii)$,  $Y_1,\dots Y_r$ are not a complete set of irreducible representations of $\G$ while, on the other hand, any irreducible representation of $\G$ occurs in the regular representation with non-zero multiplicity.  So the matrix for the system in  part $II)$ of Theorem \ref{mainth} has maximal rank and the space of solutions has dimension $\nu-r$.\epf

Combining now the results of this section we obtain Theorem \ref{mainth}.

\section{Some examples: the cyclic and  dicyclic groups}\label{example}

To clarify our results we will recover more explicit equations and make some remarks in the cyclic and  dicyclic case (which is the easiest  non-commutative case). 

\textbf{I) The cyclic case.}

Let's consider $\Gamma=\mathcal{C}_{n+1}=\{\alpha,\dots,\alpha^{n+1}=1\}$, the cyclic group of order $n+1\geq 2$. The abelian group $\mathcal{C}_{n+1}$ has, of course, $n+1$ conjugacy classes (one for each element $\alpha^s$), and $n+1$ distinct irreducible representations, all of dimension $1$. Let  $\e$ be a primitive $n+1$th root of unity, then we can write the character table of $\mathcal{C}_{n+1}$ as follows: 

\begin{center}
\begin{tabular}{c|c|c|c|c|c|}  
& $1$    & $\alpha$ & $\alpha^2$ & \dots & $\alpha^n$ \\ \hline
$\chi_1$ & $1$    & $\e$     & $\e^2$     & \dots & $\e^n$       \\ \hline
$\chi_2$ & $1$    & $\e^2$   & $\e^4$     &\dots  & $\e^{2n}$     \\ \hline
\dots & \dots  & \dots    & \dots      & \dots & \dots        \\ \hline
$\chi_n$ & $1$ &    $\e^n$   & $e^{2n}$     & \dots  & $e^{n^2}$     \\ \hline
$\chi_{n+1}$ & $1$    &    $1$   & $1$        &\dots  & $1$     \\ \hline
\end{tabular}
\end{center}

 We will write  $\chi=\chi_{1}^{\otimes N_1}\otimes\dots\otimes \chi_{n+1}^{\otimes N_{n+1}}$ for an irreducible representation of $\mathcal{C}_{n+1}^N$. $c_s$  will denote the value of the class function $c$ on $\alpha^s$.  For any fixed $N$ we can now give the  list of irreducible representations of $S_N\ltimes \mathcal{C}_{n+1}^N$ that can be extended to representations of $H_{1,k,c}(S_N\ltimes \mathcal{C}_{n+1}^N)$ for some values of $k\neq 0$:
\begin{enumerate} 
\item representations $W\otimes {\chi_{h}}^{\otimes N}$ of trivial type $(0,\dots,\stackrel{h}{N},\dots 0)$, where  $W$ has rectangular Young diagram  of size $a\times b$, can be extended for values of $(k,c)$ belonging to the hyperplane:
$$
1+k\frac{(b-a)(n+1)}{2}+\sum_{s=1}^{n}c_s \epsilon^{hs}=0\,;
$$
\item looking at Figure \ref{cyclicfig}, we observe that  representations  $W\otimes \chi\uparrow$ of non-trivial type $(N_1,\dots, N_{n+1})$ occur if and only if  $r=\#\{N_h \mbox{\ s.t.\ } N_h\neq 0\}\leq [\frac{N+1}{2}]$, where $[-]$ is the \emph{Integer Part} function (there must be $r$ non-adjacent vertices on the diagram $\tilde{A}_n$). In this case we must have $W=W_1\otimes\dots\otimes W_{n+1}$ where $W_h$ has rectangular young diagram of size $a_h\times b_h$, $\forall h$. Such representations  extend for values of the parameter $c$ belonging to the intersection of the hyperplanes:
$$
1+k\frac{(b_h-a_h)(n+1) }{2}+\sum_{s=1}^{n}c_s \epsilon^{sh}=0\,,\qquad \forall h \mbox{\ s.t.\ } N_{h}\neq 0\, .
$$
\end{enumerate}

\textbf{II) The dicyclic case.}

Let $\G=\mathcal{D}_{n-2}$, $n\geq 4$, the dicyclic group of order $4(n-2)$. This group has the following presentation:
$$
\mathcal{D}_{n-2}=\{\alpha,\beta|\;\alpha^{2(n-2)}=\beta^4=e,\quad \beta^2=\alpha^{n-2},\quad \beta\,\alpha =\alpha^{-1}\,\beta\}\, .
$$

We recall that the $n+1$ conjugacy classes of $\mathcal{D}_{n-2}$ are the following:
\begin{itemize}
\item each element of type $\alpha^s$ forms a conjugacy class together with its inverse $\alpha^{-s}$ for a total of $n-1$ classes ($e$ and $\alpha^{n-2}$ form distinct conjugacy classes by themselves);\\
\item elements of type $\beta\alpha^s$ form two distinct conjugacy classes of order $(n-2)$  for $s$ odd and $s$ even respectively.
\end{itemize}

So $\mathcal{D}_{n-2}$ has $n+1$ irreducible representation, $4$ of dimension $1$ and $n-3$  of dimension $2$. The representations of dimension $2$ are the same in the case $n-2$ even  and  $n-2$ odd  and are obtained by induction from the irreducible (non self-conjugate) representation of the normal subgroup $\mathcal{C}_{2(n-2)}$. With a proper choice of basis they  are given by:
$$
\alpha\longrightarrow \left(\begin{array}{cc}\epsilon^{h}& 0\\
                                                      0 & \epsilon^{-h}
                                     \end{array}\right),\quad
\beta\longrightarrow \left(\begin{array}{cc}0 & \e^{{h}(n-2)}\\
                                                      1 & 0
                                      \end{array}\right)
$$
for $h \in[1,n-3]\;$, and $\e$ a root of unity of order $2(n-2)$.
Some of the four  $1$-dimensional  representations, instead, differ in the case $n-2$ even and $n-2$ odd. 

If we denote by $\e$ a primitive $2(n-2)$th root of unity and by $i$ the imaginary unit  we have the following character tables for $\mathcal{D}_{n-2}$.

When  $n-2$ is even:

\begin{center}
\begin{tabular}{c|c|c|c|c|c|c|c|}
 & $1$    & $\alpha^{\pm 1}$ & $\alpha^{\pm 2}$ & \dots & $\alpha^{n-2}$ & $\beta\alpha^{2s}$ & $\beta\alpha^{2s+1}$ \\
\hline    
 $|C_h|$   & $1$    & $2$ & $2$ & \dots & $1$ & $n-2$ & $n-2$ \\ \hline      
$\begin{array}{c}\delta_h\\
1\leq h \leq n-3 \end{array}$ & $2$ & $\e^h+\e^{-h}$ & $\e^{2h}+\e^{-2h}$ & \dots & $(-1)^h2$ & $0$ & $0$ \\ \hline
$\delta_{n-2}$ & $1$ & $1$ & $1$  &\dots  & $1$ & $1$ & $1$     \\ \hline
$\delta_{n-1}$ & $1$ & $1$ & $1$  & \dots & $1$ & $-1$ & $-1$ \\ \hline
$\delta_{n}$ & $1$ & $-1$ & $1$ & \dots & $1$  & $-1$ & $1$ \\ \hline
$\delta_{n+1}$ & $1$ & $-1$ & $1$ & \dots & $1$ & $1$  & $-1$    \\ \hline
\end{tabular}
\end{center}

When  $n-2$  is odd:

\begin{center}
\begin{tabular}{c|c|c|c|c|c|c|c|}
      & $1$    & $\alpha^{\pm 1}$ & $\alpha^{\pm 2}$ & \dots & $\alpha^{n-2}$ & $\beta\alpha^{2s}$ & $\beta\alpha^{2s+1}$ \\ 
\hline    
 $|C_h|$   & $1$    & $2$ & $2$ & \dots & $1$ & $n-2$ & $n-2$ \\ 
\hline
$\begin{array}{c}\delta_h \\
1\leq h \leq n-3 \end{array}$ & $2$ & $\e^h+\e^{-h}$ & $\e^{2h}+\e^{-2h}$ & \dots & $(-1)^{h}2$ & $0$ & $0$ \\ \hline
$\delta_{n-2}$ & $1$ & $1$ & $1$  &\dots  & $1$ & $1$ & $1$     \\ \hline
$\delta_{n-1}$ & $1$ & $1$ & $1$  & \dots & $1$ & $-1$ & $-1$ \\ \hline
$\delta_{n}$ & $1$ & $-1$ & $1$ & \dots & $-1$  & $-i$ & $i$ \\ \hline
$\delta_{n+1}$ & $1$ & $-1$ & $1$ & \dots & $-1$ & $i$  & $-i$    \\ \hline
\end{tabular}
\end{center}

 We will keep  the  notation $\delta_h$ for the irreducible representation of $\mathcal{D}_{n-2}$ corresponding to this character, and we will write $\delta$ for an irreducible representation of $\mathcal{D}_{n-2}^N$ . We will  denote the values of $c$ in the following way: $c_{s}=c(\alpha^s)=c(\alpha^{-s})$ (for $s=1,\dots, (n-2)$), $c_{o}=c(\beta\alpha^{2s+1})$, $c_{e}=c(\beta\alpha^{2s})$.
We can now write the list of irreducible representations of $S_N\ltimes \mathcal{D}_{n-2}^N$ that extend to representations of the corresponding symplectic  reflection algebra for $k\neq 0$:

\begin{enumerate}
\item representations of trivial type $W\otimes \delta_h^{\otimes N} $, where $W$ has rectangular Young diagram of size $a\times b$,  extend for values of the parameter $(k,c)$ satisfying the corresponding equation among the following ones:
\begin{itemize}
\item for $\delta_h$, $ 1\leq h\leq n-3$
\begin{equation}\label{maindi1}   
1+k\,(b-a)\,(n-2)+\sum_{s=1}^{n-3}c_{s}(\e^{sh}+\e^{-sh})+c_{n-2}(-1)^h=0\nonumber\,;
\end{equation}
\item for $\delta_{n-2}$,  $\delta_{n-1}$
\begin{equation}\label{maindi2}   
 1+2k\,(b-a)\,(n-2)+2\sum_{s=1}^{n-3}c_{s}+(-1)^{n-2}c_{n-2}\pm\,(n-2)\,(c_{o}+c_{e})=0\nonumber
\end{equation}
with $+$ in the case $\delta_{n-2}$,  and $-$ in the case $\delta_{n-1}$;
\item for  $\delta_{n}$,$\delta_{n+1}$ and  $n-2$ even
\begin{equation}\label{maindi4}   
1+2k\,(b-a)\,(n-2)+2\sum_{s=1}^{n-3}(-1)^sc_{s}+c_{n-2}\mp \,(n-2)\,(c_{e}-c_{o})=0\nonumber
\end{equation}
with $-$ in the case $\delta_{n}$ and $+$ in the case $\delta_{n+1}$;
\item for $\delta_n$,  $\delta_{n+1}$ and $n-2$ odd
\begin{equation}\label{maindi3}   
1+2k\,(b-a)\,(n-2)+2\sum_{s=1}^{n-3}(-1)^s c_{s}-c_{n-2}\mp i\,(n-2)\,(c_{e}-c_{o})=0\nonumber
\end{equation}
with $-$ in the case $\delta_n$ and $+$ in the case $\delta_{n+1}$;
\end{itemize}

\item looking at Figure \ref{dicyclicfig}, we can see that representations $W\otimes \delta\uparrow$ of non-trivial type $(N_1,\dots,N_{n+1})$ occur if and only if $r=\#\{N_h \mbox{\ s.t.\ } N_h\neq 0\}\leq 4+[\frac{n-4}{2}]$. In this case we must have  $W=W_1\otimes\dots\otimes W_{n+1}$ where $W_h$ has rectangular Young diagram of size $a_h\times b_h$. These representations can be extended for values of $(k,c)$ satisfying the corresponding  above equations for any irreducible representation $\delta_h$  occurring in the type of $\delta$.
\end{enumerate}

{\bf Acknowledgments.} I am very grateful to Pavel Etingof for introducing me to the subject and for many useful conversations and comments.

\end{document}